# Centre for Mathematical Sciences India (CMS):

# Professor A.M. Mathai's 75th Birthday

By Hans J. Haubold, United Nations, CMS Member since 1983

Abstract. A brief overview on the Centre for Mathematical Sciences India, established in 1977, and its teaching and research programme is given.

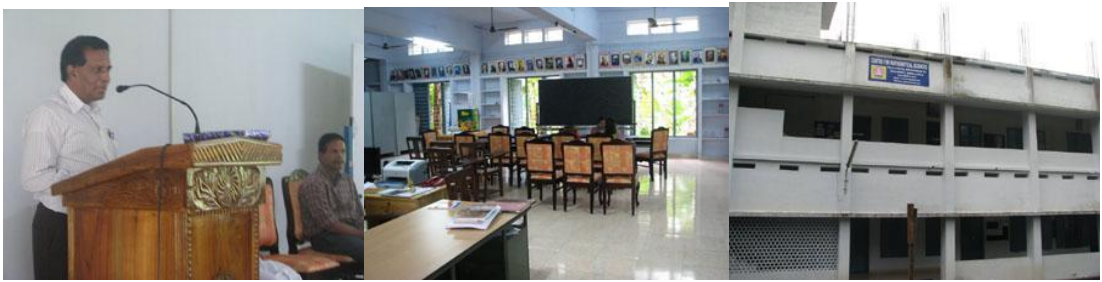

The Centre for Mathematical Sciences (CMS) was established in 1977 and registered in Trivandrum, Kerala, India, as a non-profit scientific society and a research and training centre covering all aspects of mathematics, statistics, mathematical physics, computer and information sciences. Since 1977 CMS had executed a large number of research and training projects for various central and state governmental agencies.

CMS has a publications series (books, proceedings, collections of research papers, lecture notes), a newsletter of two issues per year, a mathematics modules series (self-study books on basic topics)

Module 6: Basic Probability and Statistics: Part 1 Probability and Random Variables, December 2009, CMS, pp. 315.
Module 5: Integrals and Integration, April 2008, CMS, pp. 146.
Module 4: Limits, Continuity, Convergence and Differential Calculus, June 2008, CMS, pp. 154.
Module 3: Linear Algebra: Part III Applications of Matrices and Determinants, January 2009, CMS, pp. 507.
Module 2: Linear Algebra: Part II Determinants and Eigenvalues, June 2008, CMS, pp. 356.
Module 1: Linear Algebra: Part 1 Vectors and Matrices, March 2008, CMS, pp. 185.

and a mathematical sciences for the general public series (see cover photos below of some of the series issues). The latest books from CMS are: A.M. Mathai and H.J. Haubold (2008), *Special Functions for Applied Scientists*, Springer, New York, and A.M. Mathai, R.K. Saxena, and H.J. Haubold (2010), *The H-Function: Theory and Applications*, Springer, New York..

In 2002, CMS Pala Campus was established in a two floors finished building donated to CMS by the Diocese of Palai in Kerala, India. In 2006, Hill Area Campus of CMS was established. The office, CMS library, and most of the facilities are at CMS Pala Campus. The other campuses, namely, the South Campus (or Trivandrum campus) and the Hill Area Campus have occasional activities and libraries are being developed at South and Hill Area campuses also.

Starting from 1985, Professor Dr. A.M. Mathai of McGill University, Canada, is the Director of CMS. After taking early retirement in 2000, Professor Mathai is spending most of time at CMS and directing various CMS activities in an honorary capacity. CMS library is being built up by using the books and journals donated by Professor Mathai's colleagues, friends, and well-wishers in Canada and USA. CMS has the best library in Kerala, India, in mathematical sciences.

By the end of 2006 the Department of Science and Technology, Government of India (DST) gave a development grant to CMS. Thus, starting from December 2006 CMS is being developed as a DST Centre for Mathematical Sciences. DST has similar centres at three other locations in India.

From 1977 to 2010, CMS activities are carried out by a group of researchers in Kerala, mostly retired professors, through voluntary service. Starting from 2007 DST created full time salaried positions of three Assistant Professors, one Full Professor and one Liaison Officer. They are residing at CMS Pala Campus. DST approved up to 17 junior and senior research fellows (JRF/SRF). They are Ph.D students at CMS Pala Campus. They will receive their Ph.D degrees from Mahatma Gandhi University (MG University or Banaras Hindu University or Anna University Coimbatore), after fulfilling the residence requirements. All the JRFs and SRF are publishing papers accepted/published in international refereed journals.

CMS conducts a five-week research orientation course, called SERC School, every year. The main theme for the first sequence of five schools was special functions and functions of matrix argument and their applications.

SERC 1: Lecture Notes not published.
SERC 2: Lecture Notes Special Functions and Functions of Matrix Arguments, July 2000, CMS Publication No. 31, pp. 309.
SERC 3: Lecture Notes Special Functions and Functions of Matrix Argument: Recent Developments and Recent Applications in Statistics and Astrophysics, February 2005, CMS Publication No. 32, pp. 262.
SERC 4: Lecture Notes Special Functions and Functions of Matrix Argument: Recent Advances and Applications in Stochastic Processes, Statistics and Astrophysics, February 2006, CMS Publication No. 33, pp. 325.
SERC 5: Lecture Notes Special Functions and Functions of Matrix Argument: Recent Advances and Applications in Stochastic Processes, Statistics, Wavelet Analysis and Astrophysics, February 2007, CMS Publication No. 34, pp. 358.

The theme for the second sequence of five schools is multivariable and matrix variable calculus, statistical distributions and model building.

SERC 6: Lecture Notes Matrix Variable Calculus and Statistical Distribution Theory and Applications in Data Analysis, Model Building and Astrophysics Problems, April 2008, CMS Publication No. 36, pp. 181.

SERC 7: Lecture Notes Model Building: Multivariable and Matrix Variable Calculus with Applications Including Astrophysics, March 2009, CMS Publication No. 38, pp. 191.

SERC 8: Lecture Notes Multivariable and Matrix Variable Calculus and Applications: Stochastic Models, February 2010, CMS Publication No. 40, pp. 186.

The third school in the second sequence was held in 2010 at CMS Pala Campus. The total number of seats in each School is 30. International participation is encouraged on the expense of such foreign scholars. Local hospitality and study materials are provided free of charge by CMS. All expenses are paid by DST for the nationally selected 30 participants in each School.

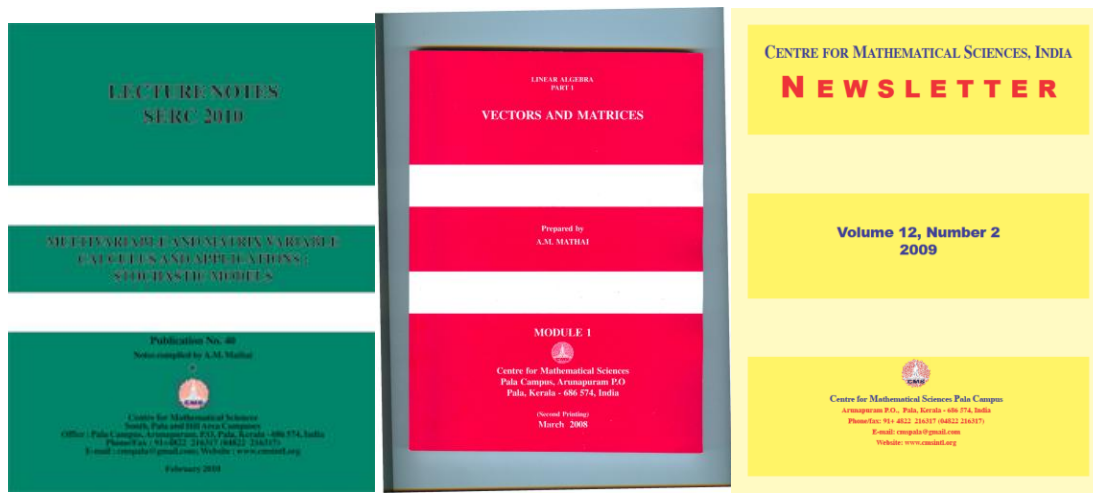

CMS is offering also an activity at the undergraduate level. There are four courses in each year covering all topics of undergraduate mathematics. Thirty students are selected by CMS from among the names recommended by the college Principals in Kerala. Each course is of 10-days duration of around 40 hours of lectures and 40 hours of problem-solving sessions. All expenses of the selected participants are met by DST.

Apart from the above two regular activities, there are lecture series of 3 days and 6 days duration by international visiting faculty.

The Director of CMS has acted as advisor in mathematics and statistics to the series of annual workshops on basic space science, organized by the United Nations, European Space Agency, National Aeronautics and Space Administration of the United States, and Space Exploration Agency of Japan, hosted by countries around the world, since their inception in 1991.

Visitors' program

There are three categories of CMS visitors: Distinguished international visitors, distinguished national visitors, faculty and students from other institutions, colleges and universities in India. For all visitors, local hospitality is provided by CMS. Those who wish to visit CMS Pala Campus need to write to the Director of CMS, giving the approximate dates and time that they would like to visit. The Director will then issue a formal invitation as per the availability of local accommodation. The recent, international visitors include Dr. A.A. Kilbas of Belarussian University, Belarus, Dr. Hans J. Haubold of the Office for Outer Space Affairs of the United Nations, Vienna, Austria, Dr. Serge B. Provost of the University of Western Ontario, Canada, Dr. Peter Moschopoulos of the University of Texas at El Paso, USA, Dr. Allan Pinkus of Israel, Dr. Francesco Mainardi of Italy, and Dr. R. Gorenflo of Germany.

Apart from the academic atmosphere at CMS, the visitors will enjoy the natural beauty all around. Pala (Palai) is situated in the heartland of agricultural activities and it is the epicenter of the spices growing region of India. Palai area grows black/white pepper, ginger, turmeric, nutmeg, cardamom, clove, coco, vanilla, cinnamon etc besides the cash crops such as coffee, tea, rubber, coconut and arecanut. The spices from here attracted the Arabs and then the Europeans to India.

The following are the contact numbers and address of CMS:
E-mail [ mathai@math.mcgill.ca; cmspala@gmail.com];
Website: www.cmsintl.org ;
Phone/fax: 91+ 4822 216317 (04822 216 317 from within India);
Postal address: Centre for Mathematical Sciences Pala Campus, Arunapuram P.O., Pala, Kerala-686574, India.